\documentclass[12pt]{amsart}
\usepackage{amssymb}
\usepackage[dvips]{graphics}
\ifx\mathrm\undefined\let\mathrm\rm\fi
\ifx\mathbf\undefined\let\mathbf\bf\fi
\ifx\mathfrak\undefined\let\mathfrak\frak\fi
\ifx\mathcal\undefined\let\mathcal\cal\fi
\ifx\mathbb\undefined\let\mathbb\Bbb\fi
\ifx\emph\undefined\let\emph\it\fi
 at9.98pt

\newcommand{\Z}{{\mathbb Z}}

\newcommand{\R}{{\mathbb R}}
\newcommand{\C}{{\mathbb C}}

\newcommand{\Ref}[1]{{(\ref{#1})}}

\newcommand{\la}{\lambda}
\newcommand{\ka}{\kappa}

\newcommand{\dontprint}[1]
{\relax}

\newtheorem%
{thm}{Theorem}[section]
\newtheorem%
{proposition}[thm]{Proposition}
\newtheorem%
{lemma}[thm]{Lemma}
\newtheorem%
{lemmadef}[thm]{Lemma-Definition}
\newtheorem%
{corollary}[thm]{Corollary}
\newtheorem%
{conjecture}[thm]{Conjecture}

\newcommand{\bea}{\begin{eqnarray*}}
\newcommand{\eea}{\end{eqnarray*}}
\newcommand{\bean}{\begin{eqnarray}}
\newcommand{\eean}{\end{eqnarray}}

\newcommand{\nc}{\newcommand}
\nc{\on}{\operatorname}
\nc{\al}{\alpha}
\nc{\ri}{\rangle}
\nc{\lef}{\langle}
\nc{\W}{{\mathcal W}}
\nc{\La}{\Lambda}
\nc{\ep}{\epsilon}
\nc{\Om}{\Omega}
\newcommand{\be}{\begin{displaymath}}
\newcommand{\ee}{\end{displaymath}}

\nc{\PCr}{{ \Bbb P  (\C[x])^r   }}

\def\ee{\emptyset}

\def\a{{\alpha}}
\def\b{{\beta}}
\def\d{{\delta}}

\newcommand{\Gm}{\Gamma}
\newcommand{\gm}{\gamma}

\newcommand{\9}{{_\infty}}

\newcommand{\8}{{\infty}}
\newcommand{\tht}{{\theta}}

\begin{document}

\title[Selberg Integrals]
{Selberg Integrals}

\author[{}]{
Alexander Varchenko ${}^{*,1}$}

\thanks{${}^1$ Supported in part by NSF grant DMS-0244579}

\begin{abstract}
The paper is written for Kluwer's Encyclopaedia of Mathematics.

\end{abstract}

\maketitle 
\centerline{\it ${}^{*}$Department of Mathematics, University of
  North Carolina at Chapel Hill,} \centerline{\it Chapel Hill, NC
  27599-3250, USA} \medskip

\centerline{August, 2004}

\bigskip

${}$
\bigskip

%\noindent
In 1944 A. Selberg proved the formula \cite{S} :
\bean\label{Selberg}
&
\int _{ \Delta^k [0,1]}\
\prod_{a=1}^k  t_a^{\alpha - 1} (1 - t_a)^{\beta -1} \prod_{1\leq a<b \leq k }
(t_a-t_b)^{2\gamma}\, dt_1 \dots dt_k 
= 
\phantom{aaaaaaaaaa}
\\
&
\phantom{aaaaaaaaaaaaaaaaaaaaaaaa}
 \prod_{j=0}^{k-1}\ 
\frac{\Gamma((j+1)\gamma)}{\Gamma(\gamma)}\
{\Gamma(\al+j \gamma)\,\Gamma(\beta+j\gamma)\, 
\over\Gamma\bigl(\alpha+\beta+(2k-2-j)\gamma\bigr)}\ ,
\notag
\eean
where $\Delta^k[x,y] = \{ t \in \R^k \ \vert\ x \leq t_k \leq \dots \leq t_1 \leq y\}$,
Re $\alpha > 0$,  Re $\beta > 0$, \linebreak
 Re $\gamma >  {\rm min}\,\{
1/k, ({\rm Re}\,\alpha)/(k-1), ({\rm Re}\,\alpha)/(k-1)\}$.
This  is a generalization of Euler's formula, 
\bea
&
\int_0^1  t^{\alpha-1}(1-t)^{\beta-1}\ dt\ =\ \frac 
{\Gamma(\alpha)\Gamma(\beta)}{\Gamma(\alpha + \beta)} \ .
\eea
The Selberg integral is one of most remarkable multi-dimensional
hypergeometric functions with many applications, see  references. 

\bigskip

Taking a suitable limit of the integral one gets
the exponential Selberg integral \cite{M}:
\bea\label{expo}
&
\phantom{a}
\int_{\Delta^k[0, \infty]}
\prod_{a=1}^k \,e^{-t_a} t_a^{\alpha-1}\!\!\prod_{1\leq a<b \leq k }
(t_a - t_b)^{2\gamma}dt_1 \dots dt_k\
= \
\prod_{j=0}^{k-1}\,{\Gamma((j+1)\gamma)\,\Gamma(\alpha+j\gamma)
\over\Gamma(\gamma)}\  .
\eea

\bigskip

There are the
Mellin-Barnes type  Selberg integrals \cite{G, TV1}:
\bea
&
\int_{-i\infty}^{i\infty}\!\!\!\cdots\!\!\int_{-i\infty}^{i\infty}
\prod_{a=1}^k \Bigl(\Gm(\a+t_a)\Gm(\b+t_a) \Gm(\gamma-t_a)\Gm(\delta-t_a)
\ \times
\\
&
\prod_{b=1}^{a-1}\,{\Gm(t_b-t_a+\epsilon)\Gm(t_a-t_b+\epsilon)
\over\Gm(t_b-t_a)\Gm(t_a-t_b)} \Bigr)\,dt_1\dots dt_k =
\\
&
= (2\pi i)^k\ k!\  \prod_{j=0}^{k-1} \Bigl({\Gm((j+1)\epsilon)\over\Gm(\epsilon)}\ 
 {\Gm(\a+\gm+j\epsilon)\Gm (\a+\d+j\epsilon)
\Gm(\b+\gm+j\epsilon)\Gm(\b+\d+j\epsilon)\over\Gm(\a+\b+\gm+\d+(2k-2-j)\epsilon)}
\Bigr) \ ,
\eea
\bea
&
\int_{-i\infty}^{i\infty}\!\!\!\cdots\!\!\int_{-i\infty}^{i\infty}\; 
\prod_{a=1}^k\,\Bigl(u^{2t_a} \Gm(\a+t_a) \Gm(\a-t_a)\,
\prod_{b=1}^{a-1}\,{\Gm(t_b-t_a+\gm)\Gm(t_a-t_b+\gm)\over\Gm(t_b-t_a)\Gm(t_a-t_b)}
\Bigr)\,dt_1 \dots dt_k  =
\\
&
(2\pi i)^k\ k!\ (u+u^{-1})^{-k (2\a+(k-1)\gm)}\,
\prod_{j=0}^{k-1}\,{\Gm((j+1)\gm)\, \Gamma (2\a+j\gm)\over\Gm(\gm)}\ ,
\eea
where $\rm{Re}\, \a,\b,\gm,\d,\epsilon, u >0$.
In particular,
\bea
\int_{-i\infty}^{\,i\infty}\Gm(\a+t)\Gm(\b+t)\Gm(\gm-t)\Gm(\d-t)\ dt  & \ =\
2\pi i \ {\Gm(\a+\gm) \Gm(\a+\d) \Gm(\b+\gamma) \Gm(\b+\d)\over\Gm(\a+\b+\gm+\d)} ,
\eea
\bea
\int_{-i\infty}^{i\infty} \Gm(\a+t) \Gm(\a-t) \ u^{2t} \ dt &\ =\ 
2\pi i \ \Gm(2\a)\ (u+u^{-1})^{-2\a} ,
\eea
which are formulae for the classical Barnes integrals \cite{WW}.

\bigskip

Let $q\in \C$, ${0<|q|<1}$, \ $(u)\9 =
\prod_{n=0}^\8(1 - q^n u)$, \
${\tht(u)=(u)\9(qu^{-1})\9(q)\9}$, \
\linebreak
$\Bbb T^k = \{ t\in \C^k\ | \ |t_1| = 1, \dots , |t_k| = 1 \}$.
The $q$-Selberg integrals have the form \cite{TV2}:
\bea
&
\int_{\Bbb T^k}\prod_{a=1}^k
\frac{\theta( q t_a / \epsilon )\, \theta( u^{k-1} \gamma \delta \epsilon t_a )}
{t_a \,(\gamma t_a)\9 \,(\delta t_a)\9 \,(\alpha/ t_a)\9\, (\beta/ t_a)\9}
\prod_{a=1}^k\prod_{b=1,\,b\neq a}^k 
\frac{ ( t_a / t_b )\9} {( u t_a / t_b )\9}\, dt_1 \dots dt_k \ =
\\
&
(2 \pi i )^k\, k!\, 
\prod_{j=0}^{k-1}
\frac{ (u)\9 \,(u^{k+j-1}\alpha \beta \gamma \delta)\9\,
\theta(u^j \gamma \epsilon)\,\theta(u^j \delta \epsilon)}
{(u^{j+1})\9\, (u^j \alpha \gamma)\9\, (u^j \beta \gamma)\9\,
(u^j \alpha \delta)\9 \,(u^j \beta \delta)\9 \,(q)\9}
\eea
where $|\alpha|<1$, $|\beta|<1$, $|\gamma|<1$, $|\delta|<1$, $|u|<1$, \ and
\bea
&
\int_{\Bbb T^k}\prod_{a=1}^k
\frac{\theta(\epsilon t_a)}
{t_a \,(\gamma t_a)\9 \,(\delta t_a)\9 }
\prod_{a=1}^k\prod_{b=1,\,b\neq a} 
\frac{ (t_a/t_b)\9} {(ut_a/t_b)\9}\, dt_1\dots dt_k \ =
\\
&
(2 \pi i )^k\, k!\, 
\prod_{j=0}^{k-1}
\frac{ (u)\9 \,(u^{j}\gamma\delta)\9\, (qu^j\gamma/\delta)\9}
{(u^{j+1})\9\, (u^j\gamma\delta)\9}
\eea
where $|\gamma|<1, |\delta|<1, |u|< 1$. In particular, 
\bea
&
\int_C\frac{\theta (qt/\epsilon)\,\theta(\gamma\delta\epsilon t)}
{(\gamma t)\9\,(\delta t)\9\,(\alpha / t)\9\, (\beta / t)\9}\,
\frac{dt}t\ =\
2\pi i \, \frac{( \alpha \beta \gamma \delta )\9 \,\theta(\gamma\epsilon)\,
\theta(\delta\epsilon)}
{(q)\9\,(\alpha \gamma)\9\,(\beta \gamma)\9\,
\9\,(\alpha \delta)\9\,(\beta \delta)\9}\ ,
\eea
where $C$ is an anti-clockwise oriented contour around the origin
$t=0$ separating the sets $\{q^s/\gamma,\, q^s/\delta\ | \ s\in \Z_{\leq 0}\}$, \
$\{q^s\alpha,\, q^s\beta \ | \ s\in \Z_{\geq 0}\}$,
\ and
\bea
&
\int_C\frac{\theta (\epsilon t)}
{(\gamma t)\9\,(\delta/ t)\9}\,
\frac{dt}t\ =\
2\pi i \, \frac{(q\gamma/\epsilon)\9 (\delta\epsilon)\9}
{(\gamma\delta)\9} \ ,
\eea
where $C$ is an anti-clockwise oriented contour around the origin
$t=0$ separating the sets $\{q^s/\gamma\ | \ s\in \Z_{\leq 0}\}$, \
$\{q^s\delta\ | \ s\in \Z_{\geq 0}\}$, see \cite{GR}.

\bigskip

One can calculate the integrals over $\Bbb T^k$ by residues and get 
formulas for the q-Selberg Jackson integrals \cite{As, E, Ka1, Ka2, TV2}.
Set 
\bea
&
A(t_1, \dots , t_k; u) \ = \ 
\prod_{a=1}^k \frac {t_a\,(q t_a/\gamma)\9\,(q t_a/\delta)\9}
{(\alpha t_a)\9\,(\beta t_a)\9}
\prod_{1\leq a < b \leq k}
\frac{(1 - t_b/t_a)\,(q u^{-1} t_b/t_a)\9}
{(u t_b / t_a)\9 }\ .
\eea
If $|qu^{k-1}| < 1$, then
\bea
&
\sum_{a=0}^k \sum_{r_1, \dots, r_k = 0}^\infty\
 (-1)^a\,
u^{\sum_{b=1}^k (k-b)(k-b+1)r_b - (k-a-1)(k-a)(1+2\sum_{b=1}^a r_b)/2}
\prod_{c=0}^{k-a-1}
\frac {\theta(u^{a+c} \gamma/\delta)}
{\theta(u^{a-c}\gamma/\delta)}
\\
&
\times\ 
A(q^{r_1}\gamma, q^{r_1+r_2}u\gamma, \dots , q^{r_1+\dots +r_a}u^{a-1}\gamma,
q^{r_{a+1}}\delta, \dots , q^{r_{a+1}+\dots +r_k}u^{k-a-1}\gamma; u)\ =
\\
&
\prod_{j=0}^{k-1}
\frac{(u)\9\, (u^{k+j-1} \alpha \beta \gamma \delta)\9\,\delta\,\theta(u^j \gamma/ \delta)}
{(u^{j+1})\9\,(u^j \alpha \gamma)\9\,(u^j \beta \gamma)\9
\,(u^j \alpha \delta)\9\,(u^j \beta \delta)\9}\ .
\eea
Set
\bea
&
B(t_1, \dots , t_k)\ =\ 
\prod_{a=1}^k
\frac{ (qt_a)\9}{(\alpha t_a)\9}
\prod_{1\leq a<b\leq k}
\frac{ (1-t_b/t_a)\,(qu^{-1}t_b/t_a)\9}
{(ut_b/t_a)\9}\ .
\eea
If $|v|< {\rm min}\, (1, |u^{k-1}|)$, then
\bea
&
\sum_{r_1, \dots, r_k = 0}^\infty\
v^{\sum_{a=1}^k (k-a+1)r_a} \
u^{-\sum_{a=1}^k (a-1)(k-a+1)r_a}\ 
B(q^{r_1},q^{r_1+r_2}u, \dots , q^{r_1+\dots +r_k}u^{k-1})\ = 
\\
&
\prod_{j=0}^{k-1}
\frac{ (u)\9\,(u^j\alpha v)\9\,(q)\9}
{(u^{j+1})\9\,(u^j\alpha)\9\,(u^{-j}v)\9} \ .
\eea

\bigskip

Let 
\bea
\theta_1(t ,\tau) =  -  
\sum_{j\in {\Z}}\ 
e^{\pi i(j+\frac{1}{2})^{2}\tau+2\pi i(j+\frac{1}{2})(t +\frac{1}{2})},
\qquad
\theta_{\ka,n}(t,\tau) = 
\sum_{j\in\mathbb{Z}}\ e^{2\pi i\ka(j+\frac{n}{2\ka})^{2}\tau+2\pi
  i\ka(j+\frac{n}{2\ka})t}
\eea
where $\kappa, n \in \Z$,  $\ka \geq 2$.  Let
\bea
&
\sigma_{\la}(t,\tau)=
\frac{\theta_1(\la-t,\tau)\ \theta'_1(0,\tau)}
{\theta_1(\la,\tau)\ \theta_1(t,\tau)},
\qquad
E(t,\tau)=\frac{\theta_1(t,\tau)}{\theta'_1(0,\tau)},
\eea
where $'$ denotes the derivative with respect to the first argument.
For $k \in \Z_{>0}$, the elliptic Selberg integral is the integral
$I_k(\lambda,\tau) = J_k(\lambda,\tau) + (-1)^{k+1}J_k(-\lambda,\tau)$, where
\bea
&
J_k(\lambda,\tau) \ = \ 
\int_{\Delta^k[0,1]} 
\theta_{2k+2,k+1}(\lambda + \frac1{k+1}\sum_{a=1}^k t_a, \, \tau)
\times
\\
&
\prod_{a=1}^k E(t_a,\tau)^{-k/(k+1)} \sigma_\lambda(t_a,\tau)
\ \prod_{1\leq a<b\leq k}
E(t_a-t_b,\tau)^{1/(k+1)} \,dt_1 \dots dt_k \ .
\eea
Then \cite{FSV1, FSV2}
\bea
&
I_k(\lambda, \tau)\ =\
c_k \, B_k(\frac 12+\frac 1{2k+2}, -\frac k{k+1},\frac 1{2k+2})\ 
\theta_1(\lambda,\tau)^{k+1} \ ,
\eea
where
$c_k = (2\pi)^{k/2} e^{\pi i k/(k+1)} e^{-\pi i(k+2)/4}
\prod_{a=1}^k(e^{-\pi i a/(k+1)} - 1)$ and $B_k(\alpha,\beta,\gamma)$
is the Selberg integral in \Ref{Selberg}.

\bigskip

All previous versions of the Selberg integral are related to representation 
theory of the Lie algebra $sl_2$ and appear as solutions to different versions 
of the KZ equations and dynamical equations, see \cite{FMTV, TV1, TV2, TV3, 
TV4, TV5, FSV1, FSV2, V}.
The following integral is an example of a Selberg type integral associated with $sl_3$.

Let $k_1, k_2$ be non-negative integers, $k_1\geq k_2$.
For any non-decreasing map $M : \{1, \dots , k_2\} \linebreak
\to
\{ 1, \dots , k_1\}$ such that
$M(b) \leq k_1 - k_2 + b$ for any $b = 1, \dots, k_2$, introduce a domain 
$\Delta^{k_1,k_2}_M [x,y]$
in $\R^{k_1+k_2}$ with coordinates $t_1, \dots , t_{k_1}, s_1, \dots , s_{k_2}$ 
defined by the inequalities:
\bea
&
x \leq t_{k_1} \leq \dots \leq t_1 \leq y\,,
\qquad 
x \leq s_{k_2} \leq \dots \leq s_1 \leq y\,,
\\
&
t_{M(b)} \leq s_b \leq t_{M(b)-1}\,, \qquad b = 1, \dots ,  k_2\,.
\eea
Here $t_0 = y$. For $\gamma \in \C$, consider the chain
$C^{k_1,k_2}_\gamma [x,y]
\,=\, \sum_M\,X_{M,\gamma}^{k_1,k_2} \Delta^{k_1,k_2}_M[x,y]$ with coefficients
\bea
&
X_{M,\gamma}^{k_1,k_2}\ = \
\prod_{b=1}^{k_2}
\frac{ {\rm sin}\,(\pi ( k_1 - k_2 - M(b) + b + 1) \gamma )}
{ {\rm sin}\,(\pi ( k_1 - k_2 + b) \gamma )}\ .
\eea
Then \cite{TV5}
\bea
&
\int_{C^{k_1,k_2}_\gamma [0,1]}
\prod_{a=1}^{k_1}t_a^{\alpha-1}(1-t_a)^{\beta_1-1}
\prod_{b=1}^{k_2}(1-s_b)^{\beta_2-1}
\prod_{a=1}^{k_1}\prod_{b=1}^{k_2}|t_a-s_b|^{-\gamma} \ \times
\\
&
\prod_{1\leq a<b\leq k_1} (t_a-t_b)^{2\gamma}
\prod_{1\leq a<b\leq k_2} (s_a-s_b)^{2\gamma}\
dt_1 \dots dt_{k_1}\,ds_1 \dots ds_{k_2}\ =
\\
&
\prod_{j=0}^{k_1-1}
\frac{ \Gamma(\alpha + j\gamma)\,\Gamma( (j+1)\gamma) }
{\Gamma(\gamma)}
\prod_{j=0}^{k_1-k_2-1}
\frac{ \Gamma(\beta_1 + j\gamma)}
{\Gamma(\alpha + \beta_1 + (2k_1-k_2-2-j)\gamma)} \ \times
\\
&
\prod_{j=0}^{k_2-1}
\frac{ \Gamma(\beta_2 + j\gamma)\,
\Gamma(\beta_1+\beta_2 + (j-1)\gamma)\,
\Gamma(1+(j-k_1)\gamma)\,\Gamma( (j+1)\gamma)}
{\Gamma(\beta_2 + 1 + (2k_2-k_1-2-j)\gamma)\,
\Gamma(\alpha + \beta_1+\beta_2 + (k_1+k_2-3-j)\gamma)\,\Gamma(\gamma)}\ .
\eea

\bigskip

In the  $k$-dimensional Euclidean space $\R^k$ consider
the reflection hyperplanes of a Coxeter group $G$. 
Let $P(t)$ be the product
of distances of the point $t = (t_1,...,t_k)$ from all reflection hyperplanes
associated with $G$. Let $N$ be the number of hyperplanes. Then \cite{Ma, Op, M}
$$
\int_{\R^k}\
e^{-\sum t_i^2/2}\ 
|P(t)|^{2\gamma}\
d t_1 \dots d t_k \ 
	=
\ 2^{-N\gamma}\ (2\pi)^{k/2}\
\prod_{j=1}^k\
\frac{\Gamma(1+\gamma d_j)}{\Gamma(1+\gamma)}\ ,
$$
where $d_j$ are the degrees of  basis polynomials in the space of
homogeneous polynomials, which are invariant with respect to $G$.
The $q$-analogs of this formula were proved in  \cite{Ch}.

\bigskip

Here is  an example of a Selberg type integral involving the 
elliptic gamma function, see \cite{Sp}. Let
\bea
&
\Gamma (t; p,q)\ =\
\prod_{j,k = 0}^\infty
\frac{1-t^{-1}q^{j+1}p^{k+1}}{1-t     q^ j   p^{k  }}
\eea
be the elliptic gamma function and $(q;q)\9 = \prod_{n=1}^\infty (1-q^n)$.
For complex parameters $u_0, \dots , u_4$, set $A = \prod_{m=0}^4 u_m$.
Assume that $|p|<1, |q|<1, |u_m|<1, |pq|< |A|$. Then
\bea
&
\frac 1{2\pi i}\int_{\Bbb T^1}
\frac 
{\prod_{m=0}^4 \Gamma(tu_m;\, p,q)\, \Gamma(t^{-1}u_m;\,p,q)}
{\Gamma(t^2;\,p,q)\, \Gamma(t^{-2};\,p,q)\,
\Gamma(tA;\,p,q)\, \Gamma(t^{-1}A;\,p,q)}
\frac{dt} t\ =
\frac{ 2 \prod_{0\leq l<m\leq 4} \Gamma(u_l u_m;\,p,q)}
{(q;q)\9\, (p;p)\9\, \prod_{m=0}^4
\Gamma(A u_m^{-1};\,p,q)} \ .
\eea

\end{document}